\documentclass[12pt]{article} 
\usepackage{amsmath}
\usepackage{amssymb} 
\usepackage{amscd}
\usepackage{amsthm}
\usepackage[utf8]{inputenc}
\usepackage{stmaryrd}
\usepackage{hyperref}
\newtheorem{theorem}{Theorem}

\newtheorem{proposition}[theorem]{Proposition}

\def\gr{\mathrm{gr}}

\def\dim{{\mbox{dim}}}

\def\ker{{\mbox{Ker}}}     
\def\Der{{\mbox{Der}}}
\def\der{{\mbox{\tiny Der}}}

\def\Lie{{\mbox{Lie}}}

\def\cala{{\cal A}}

\def\call{{\cal L}}

\def\calh{{\cal H}}

\def\calc{{\cal C}} 
\def\cale{{\cal E}}

  \def\fracf{{\mathfrak f}} 
   \def\fracg{{\mathfrak g}}

\def\bbbone{\mbox{\rm 1\hspace {-.6em} l}}

\def\Diff{{\mbox{Diff}}}
\def\Aut{{\mbox{Aut}}}

\numberwithin{equation}{section}

\begin{document}
\enlargethispage{3cm}

\thispagestyle{empty}
\begin{center}
{\bf EXCEPTIONAL QUANTUM GEOMETRY}
\end{center} 
\begin{center}
{\bf AND PARTICLE PHYSICS}
\end{center}
\vspace{1cm}

\begin{center}
 Michel DUBOIS-VIOLETTE
\footnote{Laboratoire de Physique Th\'eorique,
 CNRS, 
Universit\'e Paris-Sud,
Universit\'e Paris-Saclay,
B\^atiment 210,
F-91 405 Orsay.\\
Michel.Dubois-Violette$@$u-psud.fr}\\
\end{center}
\vspace{0,5cm}
\begin{center}
{\sl To the memory of Raymond Stora}
\end{center}
 \vspace{0,5cm}

\begin{abstract}
Based on an interpretation of the quark-lepton symmetry in terms of the unimodularity of the color group $SU(3)$ and on the existence of 3 generations, we develop an argumentation suggesting that the ``finite quantum space" corresponding to the exceptional real Jordan algebra of dimension 27 (the Euclidean Albert algebra) is relevant for the description of internal spaces in the theory of particles. In particular, the triality which corresponds to the 3 off-diagonal octonionic elements of the exceptional algebra is associated to the 3 generations of the Standard Model while the representation of the octonions as a complex 4-dimensional space $\mathbb C\oplus\mathbb C^3$ is associated to the quark-lepton symmetry, (one complex for the lepton and 3 for the corresponding quark). More generally it is is suggested that the replacement of the algebra of real functions on spacetime by the algebra of functions on spacetime with values in a finite-dimensional Euclidean Jordan algebra which plays the role of ``the algebra of real functions" on the corresponding almost classical quantum spacetime is relevant in particle physics. This leads us to study the theory of Jordan modules and to develop the differential calculus over Jordan algebras, (i.e. to introduce the appropriate notion of differential forms). We formulate the corresponding definition of connections on Jordan modules.
 
\end{abstract}

\vfill
\noindent \today\\

\newpage
\tableofcontents
\newpage

\section{Introduction}

It is well known that the Standard Model of particles is very successful but contains several inputs which ought to have explainations at a fundamental level. Some of these inputs like the existence of the Higgs particles, which seem to have recently acquired experimental confimations, appear naturally for instance in the (almost commutative) noncommutative geometric formulations \cite{mdv-ker-mad:1989a}, \cite{mdv-ker-mad:1988b}, \cite{mdv:1991}, \cite{ac-lot:1990}, \cite{cha-ac:1997}, \cite{cha-ac-mar:2007} or in the superconnection formulations \cite{roe:2000}. Some other structural inputs come almost directly from the experimental observation. For instance the quark-lepton symmetry or the existence of 3 generations belong to these latter inputs and it is one of our aims here to connect these 2 inputs and to suggest some theoretical explanations for these facts.\\

By the quark-lepton symmetry, we mean the fact that to each quark corresponds one lepton and conversely. For instance to the quark $u$ corresponds the electronic neutrino $\nu_e$, to the quark $d$ corresponds the electron $e$, etc..\\
On the theoretical side, it is worth noticing here that the cancellation of anomalies is a very strong argument for the quark-lepton symmetry.\\

The classical quark field is a spinor field with values in a complex 3-dimensional space acted by the color group $SU(3)$. Thus, forgetting the spin etc., at each point of spacetime, the internal space $E$ for a quark is a complex 3-dimensional Hilbert space (since $SU(3)\subset U(3)$) which is endowed with a complex volume that is with an antisymmetric $\mathbb C$-trilinear form (since elements of $SU(3)$ are of determinant 1). By using the Hilbertian scalar product of $E$, one can transform the volume into an antilinear antisymmetric product $\vartimes$ on $E$. Thus $E$ is equipped with a product 
\[
\vartimes:E\times E\rightarrow E
\]
and the Hilbertian scalar product
\[
\langle, \rangle : E \times E\rightarrow \mathbb C
\]
and $SU(3)$ is the group of the $\mathbb C$-linear transformations which preserve these structures. It is then natural to combine these products into an $SU(3)$-invariant product on the Hilbertian direct sum $A=\mathbb C \oplus E$ such that  $1\in \mathbb C\subset A$ is a unit denoted by $\bbbone$ for the product of $A$ and such that the norm of the product of 2 elements of $A$ is the product of the norms of these 2 elements. The product on $A$ is then not $\mathbb C$-bilinear but is only real bilinear and the corresponding real algebra is isomorphic to the algebra of octonions \cite{bae:2002}. The group $SU(3)$ identifies then with the group of $\mathbb C$-linear isomorphisms of $A$ ($\simeq
\mathbb C^4$) which preserve the product. Under this action the invariant subspace $E$ of $A=\mathbb C\bbbone \oplus E$ corresponds to the fundamental representation of $SU(3)$ while its orthogonal corresponds to the trivial representation. Thinking of $E$ as the internal quark space, it is then natural to identify the component $\mathbb C\bbbone$ of $A$ as the trivial internal space of the corresponding lepton. In this way one connects the quark-lepton symmetry to the $S$ of $SU(3)$, i.e. to the unimodularity condition.\\

It is worth noticing here that this interpretation of the quark-lepton symmetry leads directly to the original interpretation given in \cite{gun-gur:1974} of the color group $SU(3)$ as a subgroup of the automorphism group $G_2$ of the octonion algebra. Indeed it is classical (see e.g. in  \cite{ram:1976}, \cite{bae:2002}, \cite{yok:2009}) that the automorphism group of the real algebra $\mathbb O$ of octonions is the first exceptional group $G_2$ and that $SU(3)$ identifies with the subgroup of $G_2$ of automorphisms which preserve a given imaginary unit of $\mathbb O$. This is of course directly connected to the above construction which is explained in details in Section 2; the corresponding imaginary unit of $\mathbb O$ being then the imaginary unit $i$ of $\mathbb C$.\\

Since $SU(3)$ is a gauge group, that is the structure group of a bundle over spacetime, this construction means that at each point of spacetime one has a complex 4-dimensional vector space $A=\mathbb C\oplus E$ which is a real algebra (for the underlying structure of real 8-dimensional vector space) isomorphic to the octonion algebra, ($SU(3)$ being the automorphism group of this structure). In other words associated to the $SU(3)$-bundle we have a complex vector bundle $A$ of rank 4 which is also a bundle in octonion algebras. A natural question is what is the role of this local algebra structure if it is relevant ?\\

Let us observe that there are 3 generations and that the principle of triality combined with octonions leads directly to the exceptional Jordan algebra $J^8_3=\calh_3(\mathbb O)$ of octonionic hermitian 3 by 3 matrices, see e.g. in \cite{ada:1980}. This algebra was introduced in \cite{jor-vn-wig:1934} where it was pointed out that it is a quantum version of an algebra of ``real functions" on a ``finite quantum space" (see also  \cite{gun-pir-rue:1978}) as any formally real (Euclidean) finite-dimensional Jordan algebra.\\
 
 The above discussion suggests to put over each spacetime point a finite quantum space corresponding to the exceptional Jordan algebra $J^8_3$. As before this corresponds to a bundle in algebras $J^8_3$ which is associated to the above (color) $SU(3)$-bundle. In other words this suggests to replace the real algebra $C^\infty(M)$ of smooth functions on spacetime $M$ by the real Jordan algebra $J^8_3(M)$ of smooth sections of this bundle in algebras $J^8_3$. Since in this paper we assume that $M$ is the usual Minkowski space $\mathbb R^4$, this algebra $J^8_3(M)$ is isomorphic to the algebra $C^\infty(M,J^8_3)=J^8_3\otimes C^\infty(M)$ of smooth $J^8_3$-valued functions on $M$. That is, we take   $J^8_3(M)=J^8_3\otimes C^\infty(M)$ which is interpreted as the ``algebra of real functions" on the quantum space which is the product of the ``finite quantum space" corresponding to $J^8_3$ with the spacetime $M$. Since the ``quantum part" is ``finite", we refer to such a quantum space as an ``almost classical quantum space". Classical matter fields are then elements of $J^8_3(M)$-modules and gluon fields are part of connections on these modules. However none of these two items, modules and connections, are straightforward for Jordan algebras.\\
 
There is a notion of (bi)module  \cite{sch:1995}, \cite{jac:1968} for Jordan algebras (see also in \cite{mdv:2001}) which will be used here. It seems that no differential calculus has been developed so far for Jordan algebras, a step which is certainly necessary to speak of connections. We define such a differential calculus over Jordan algebras and describe several examples and their properties. A particular attention is devoted to the case of the exceptional Jordan algebra $J^8_3$. An appropriate associated notion of connection on Jordan modules is introduced.\\
 
It is our aim to develop in the furure these notions in the present context in order to get a natural description of the field content of the Standard Model of particles and of some possible generalizations.\\

The present approach is close in spirit to the approach of noncommutative geometry. However here it is not a finite noncommutative space corresponding to a finite-dimensional associative algebra but a finite quantum space corresponding to a finite-dimensional Euclidean Jordan algebra which is added to spacetime. Furthermore since the exceptional Albert algebra is involved, one cannot use directly the technics of noncommutative geometry (i.e. of associative algebra) and some care must be taken to handle such theory. Nevertheless as in the noncommutative geometric approaches, it is expected that the apparition of the Higgs fields is a consequence of the quantum nature of the almost classical spacetime (i.e. components of connections in the quantum directions).\\

Finally let us recall that occurrence in physics of octonions and exceptional structures has a long history, (see e.g.  \cite{jor-vn-wig:1934},\cite{gun-gur:1973},\cite{gun-gur:1974}, \cite{gun-pir-rue:1978}, \cite{hor-bie:1979}, \cite{dix:1994}, \cite{ram:2003}). They appear very naturally in extensions of the Standard Model. More recently, in \cite{boy-far:2014}, \cite{far-boy:2014} and \cite{far-boy:2015}, an interesting tentative approach to incorporate nonassociativity in the framework of the spectral action principle has been carried over.\\

Our notations are standard and we use  Einstein convention of summation over repeated up-down indices in the formulas.

\section{Unimodularity of $SU(3)$ and the quark-lepton symmetry}

\subsection{$SU(3)$-spaces}
In the following, an $SU(3)$-{\sl space} will be a complex 3-dimensional Hilbert space $E$ equipped with an antisymmetric complex 3-linear form $v:\Lambda^3E\rightarrow \mathbb C$ of norm $\parallel v\parallel=1$, that is such that one has
\[
\vert v(e_1,e_2,e_3)\vert=1
\]
for any orthonormal basis $(e_k)$ of $E$.\\

An $SU(3)$-{\sl basis} of $E$ will be an orthonormal basis $(e_1,e_2,e_3)$ of $E$ such that one has $v(e_1,e_2,e_3)=1$. Given such an $SU(3)$-basis $(e_1,e_2,e_3)$, one defines a bijection of the group $SU(3)$ onto the set of all $SU(3)$-basis\linebreak[4] $g\mapsto (e_1^{(g)}, e_2^{(g)}, e_3^{(g)})$ by setting
\begin{equation}
e_k^{(g)}=e_\ell g^\ell_k
\label{basis}
\end{equation}
for $k\in \{1,2,3\}$ where $g^\ell_k$ are the matrix elements of $g\in SU(3)$ (in the fundamental representation).\\

Given an $SU(3)$-basis $(e_1,e_2,e_3)$ of $E$ one defines the components $Z^k\in \mathbb C$ of a vector $Z\in E$ by
\[
Z=e_kZ^k
\]
and $E$ can then be identified with $\mathbb C^3$ with Hilbertian scalar product given by
\begin{equation}
\langle Z_1,Z_2\rangle=\sum^3_{k=1} \bar Z^k_1 Z^k_2
\label{scal}
\end{equation}
and $v$ given by
\begin{equation}
v(Z_1,Z_2,Z_3)=\varepsilon_{k\ell m}Z^k_1 Z^\ell_2 Z^m_3
\label{vol}
\end{equation}
for $Z_1, Z_2, Z_3\in E$. An element $Z\in E$ is identified with the column $[Z^k]\in \mathbb C^3$ of its components and a change of $SU(3)$-basis corresponds to the left action of an $SU(3)$-matrix on $[Z^k]$, (which leaves invariant the forms (\ref{scal}) and (\ref{vol})).

\subsection{$SU(3)$-algebra version of octonion algebra}

One can combine the 3-linear form $v$ on $E$ with the Hilbertian scalar product to define an antibilinear antisymmetric product $(Z_1,Z_2)\mapsto Z_1\vartimes Z_2$ on $E$ by
\begin{equation}
v(Z_1,Z_2,Z_3)=\langle Z_1\vartimes  Z_2,Z_3\rangle
\label{vectp}
\end{equation}
for any $Z_1,Z_2,Z_3\in E$. In the $SU(3)$-basis $(e_1,e_2,e_3)$ and the corresponding identification of $E$ with $\mathbb C^3$ this product reads
\begin{equation}
(Z_1\vartimes  Z_2)^k=\varepsilon_{k\ell m} \bar Z^\ell_1 \bar Z^m_2
\label{prv}
\end{equation}
for $k\in \{1,2,3\}$. Notice that this product on $E$ is nonassociative and is invariant by the action of $SU(3)$. Thus one has a product
\[
\vartimes:E\times E\rightarrow E, \ \ (Z,Z')\mapsto Z\vartimes  Z'
\]
and the scalar product
\[
\langle, \rangle:E\times E \rightarrow \mathbb C,(Z,Z')\mapsto \langle Z,Z'\rangle
\]
which are invariant by the action of $SU(3)$ on $E$. This representation of $SU(3)$ is the fundamental representation and $SU(3)$ is exactly  the subgroup of the complex linear group of $E$ which preserves these two products; $\mathbb C$ being equipped with the trivial representation of $SU(3)$
\[
(g,z)\mapsto gz=z
\]
for $g\in SU(3)$, $z\in \mathbb C$. It is therefore natural to combine these products into an $SU(3)$-invariant product
\[
((z,Z),(z',Z'))\mapsto (z,Z)(z',Z')
\]
on $\mathbb C\oplus E $ in such a way that
\begin{equation}
(0,Z)(0,Z')=(\alpha\langle Z,Z'\rangle,\beta Z\vartimes Z')
\label{p0}
\end{equation}
for some $\alpha,\beta\in \mathbb C$ and that $\bbbone=(1,0)$ is a unit, i.e. that one has in particular
\[
(1,0)(0,Z)=(0,Z)(1,0)=(0,Z)
\]
for any $Z\in E$ with of course $(1,0)(1,0)=(1,0)$ and more generally
\[
(z,0)(z',0)=(zz',0)
\]
for $z,z'\in \mathbb C$.\\

If $\mathbb C\oplus E$ is endowed with its natural structure of 4-dimensional Hilbert space (Hilbertian direct sum), one has
\[
\parallel (0,Z)(0,Z')\parallel^2=(\vert\alpha\vert^2-\vert\beta\vert^2)\vert\langle Z,Z'\rangle\vert^2 + \vert \beta\vert^2\parallel Z\parallel^2\parallel Z'\parallel^2
\]
so that one has
\[
\parallel (0,Z)(0,Z')\parallel^2=\parallel (0,Z)\parallel^2 \parallel(0,Z')\parallel^2
\]
by chosing
\begin{equation}
\vert \alpha\vert^2=\vert \beta\vert^2=1
\label{nor}
\end{equation}
as normalization. One has also of course
\[
\parallel (z,0)(z',0)\parallel^2=\parallel (z,0)\parallel^2\parallel(z',0)\parallel^2
\]
so it is natural to require
\[
\parallel (z,Z)(z',Z')\parallel^2=\parallel(z,Z)\parallel^2\parallel (z',Z')\parallel^2
\]
with $\parallel(z,Z)\parallel^2=\vert z\vert^2+\parallel Z\parallel^2$ and the real bilinearity of the product.\\
A solution is given by
\begin{equation}
(z,Z)(z',Z')=(zz'-\langle Z,Z'\rangle,\bar zZ'+z'Z+iZ\vartimes Z')
\label{proc}
\end{equation}
This product is probably essentially unique up to irrelevant normalizations. Setting
\begin{equation}
\overline{(z,Z)}=(\bar z, -Z) 
\label{cjc}
\end{equation}
one obtains in particular
\[
\overline{(z,Z)} (z,Z)=(z,Z)\overline{(z,Z)}=(\vert z\vert^2+\parallel Z\parallel^2)(1,0)
\]
i.e.
\begin{equation}
\overline{(z,Z)}(z,Z)=(z,Z)\overline{(z,Z)}=\parallel (z,Z)\parallel^2 \bbbone
\label{nj}
\end{equation}
which implies that as a real 8-dimensional algebra for the product (\ref{proc}),\linebreak[4] $A=\mathbb C\bbbone\oplus E$ is a normed division algebra which is (necessarily) isomorphic to the algebra $\mathbb O$ of octonions. The group $SU(3)$ is the group of automorphisms of this algebra which preserves the structure of 4-dimensional complex vector space of $\mathbb C\oplus E$ or, equivalently, the group of all complex automorphisms of $\mathbb C\oplus E$ which preserve the product given by (\ref{proc}).

\subsection{The quark-lepton symmetry}
 
 The representation of $SU(3)$ on $A=\mathbb C\bbbone \oplus E$ is the direct sum of the fundamental representation on $E$ with the trivial representation on its orthogonal complement $\mathbb C\bbbone \simeq \mathbb C$. If one interprets $E$ as the internal space for a quark, it is therefore natural to identify $\mathbb C\bbbone$ with the trivial internal space of the corresponding lepton in the quark-lepton symmetry. This identification connects the quark-lepton symmetry to the unimodularity of the color group $SU(3)$.\\
 
 This raises the following question : Is there a reason for the occurrence of the division algebra $A\simeq \mathbb O$ in the present context ?
 
 \subsection{Gauge theoretical aspect}
 
 Up to now, the discussion is over each spacetime point. One knows however that the color group $SU(3)$ is the structure group of an $SU(3)$-gauge theory of strong interactions. Indeed the gluon field which induces the strong interactions is an $SU(3)$-connection.\\
 
 This means that one has a complex vector bundle over spacetime with fiber $\mathbb C\oplus \mathbb C^3$ associated with a principal $SU(3)$-bundle which
 is the direct sum of a trivial complex vector bundle of rank 1 with a complex vector bundle of rank 3 corresponding to the fundamental representation of $SU(3)$. This bundle is also an associated bundle in algebras with product given by (\ref{proc}) in each fiber. The gluon field being then a connection on this bundle.\\
 
 Let $A(M)$ be the algebra of smooth sections of this bundle and let us denote by $\Der(\cala)$ the Lie algebra of derivations of an algebra $\cala$. The covariant derivative $\nabla_X$ along the vector field $X\in\Gamma(TM)=\Der(C^\infty(M)$) associated to the gluon field $SU(3)$-connection is a derivation of $A(M)$ and the restriction $\rho$ of this derivation to $C^\infty(M)=C^\infty(M)\bbbone \subset A(M)$ satisfies 
 \[
 \rho(\nabla_X)=X
 \]
 which means that $X\mapsto \nabla_X=s_\nabla(X)$ is a splitting
 \[
 s_\nabla:\Der(C^\infty(M)) \rightarrow \Der(A(M))
 \]
 of the exact sequence of $C^\infty(M)$-modules (compare with \cite{mdv-mas:1998})
 \begin{equation}
 0\rightarrow \text{Int}(A(M))\stackrel{\subset}{\rightarrow} \Der(A(M))\stackrel{\rho}{\rightarrow}\Der(C^\infty(M))\rightarrow 0
 \label{seqD}
 \end{equation}
 where $\text{Int}(A(M))=C^\infty(M,\fracg_2)$ with the exceptional Lie algebra $\fracg_2$ identified with the Lie algebra $\Der(\mathbb O)$ of derivations of the octonion algebra which are all inner (see e.g. in \cite{sch:1995}) and thus $\Der(\mathbb O)= \text{Int}(\mathbb O)=\fracg_2$ and one has $\text{Int}(A(M))=C^\infty(M,\fracg_2)$. This splitting of (\ref{seqD}) gives a splitting of the exact sequence of $C^\infty(M)$-modules
 \begin{equation}
 0\rightarrow C^\infty(M,{\frak s\frak u}(3))\stackrel{\subset}{\rightarrow} \Der_0(A(M))\stackrel{\rho}{\rightarrow}\Der(C^\infty(M))\rightarrow 0
 \label{seqDo}
 \end{equation}
 where $\Der_0(A(M))$ is the $C^\infty(M)$-submodule and Lie subalgebra of the derivation of $A(M)$ compatible with the representation $\mathbb C\oplus \mathbb C^3$ of $\mathbb O$. Notice that (\ref{seqD}) and (\ref{seqDo}) are also exact sequences of Lie algebras but that the obstruction for $s_\nabla$ to be a Lie algebra homomorphism is the curvature of $\nabla$, that is the field strenght of the gluon field. Corresponding to $s_\nabla$ one defines a projection $\pi_\nabla$ of $\Der(A(M))$ onto $\text{Int}(A(M))=C^\infty(M,\fracg_2)=\ker(\rho)$ by setting (as in \cite{mdv-mas:1998})
 \begin{equation}
 \pi_\nabla(\delta)=\delta-\nabla_{\rho(\delta)}
 \label{pidelta}
 \end{equation}
 which restricts as a projection of $\Der_0(A(M))$ onto $C^\infty(M,{\frak s\frak u}(3))$.\\

 In the following it is assumed that the spacetime is the ordinary Minkowski space $\mathbb R^4$. This simplifying assumption implies that all the bundles over spacetime can be considered as products so their sections identify with functions with values in fixed spaces (like $\mathbb C\oplus\mathbb C^3$, etc.). However the whole discussion extends easily to the case where spacetime is a non-trivial pseudo-Riemannian manifold and where the bundles are non-trivial.

\section{The 3-generations, triality and the exceptional finite quantum space}

\subsection{the 3-generations}
 
 There are 6 flavours of quark-lepton 
 \[
 (u,\nu_e),\> (d,e),\> (c,\nu_\mu),\> (s,\mu),\> (t,\nu_\tau),\> (b,\tau)
 \]
  which can be grouped in 3 generations of doublets of quark-lepton which are the columns of the following table
 \begin{center}
 \begin{tabular}{|c|c|c|c|c|}\hline 
 generations & & & \\ \hline
 quarks $Q=2/3$	& \phantom{abc}$u$ \phantom{abc}& \phantom{abc} $c$ \phantom{abc}& \phantom{abc}$t$\phantom{abc}\\
 leptons $Q=0$ & $\nu_e$ & $\nu_\mu$ & $\nu_\tau$ \\  \hline
 quarks $Q=-1/3$ & $d$ & $s$ & $b$\\  
 leptons $Q=-1$ & $e$ & $\mu$ & $\tau$\\ \hline
 \end{tabular}
 \end{center}
 where $Q$ denotes the electric charge.\\
 
 This is the present experimental situation which reveals a sort of ``triality". This triality combined with the above interpretation of the quark-lepton symmetry is the starting point for the following analysis which suggests to add over each spacetime point an exceptional ``finite quantum space" corresponding to the exceptional Jordan algebra $J^8_3$ of the hermitian $3\times 3$ octonionic matrices and to take the internal spaces of the basic fermions as elements of appropriate modules over this algebra.\\
 
  In fact the combination of octonions with triality leads naturally to the exceptional Jordan algebra $J^8_3$ (see e.g. in \cite{ada:1980}, \cite{jac:1968}, \cite{yok:2009}) and it turns out that this Jordan algebra can be considered as the ``algebra of real functions" on a ``finite quantum space".\\

 Let us first make precise what is meant here by a finite quantum space.
 
 \subsection{Finite quantum spaces}
 
The notion of finite-dimensional Euclidean (or formally real) Jordan algebra was introduced and analyzed in \cite{jor-vn-wig:1934} under the name $r$-number algebra in order to formalize the properties of the observables of finite quantum systems. A minimal set of requirement for the observables of a finite quantum system is that it is a finite-dimensional real vector space $J$ such that one can consistently define powers of any element and that if a sum of squares of elements vanishes it implies that each of these elements vanishes. Given the square $x\mapsto x^2$, one defines a symmetric bilinear product on $J$  by setting $x\cdot y=1/2((x+y)^2-x^2-y^2)$ for $x,y\in J$.

 So let $J$ be a finite-dimensional real vector space endowed with a symmetric bilinear product
\begin{equation}
(x,y) \mapsto x.y = y.x
\label{prodcom}
\end{equation}
for $x,y\in J$ and assume that the following condition is satisfied
\begin{equation}
\sum_{i\in I} (x_i)^2=0 \Rightarrow x_i=0,\>\> \forall i\in I
\label{formre}
\end{equation}
for any finite family $(x_i)_{\in I}$ in $J$, where $(x)^2=x.x$ for $x\in J$. The above condition for the product is refered to as {\sl formal reality condition}. One has the following result \cite{jor-vn-wig:1934}.

\begin{theorem} \label{POJO}
Let $J$ be as above and define $x^n\in J$ for $x\in J$ and $n\in \mathbb N^+$ by induction on $n$ as $x^1=x$ and $x^{n+1}=x.x^n$.\\
 Then the following conditions $\mathrm{(i)}$ and $\mathrm{(ii)}$ are equivalent :\\
$\mathrm{(i)}\>\>\>  x^r.x^s=x^{r+s},\>\>  \text{for any}\>\> x\in J$ and $r,s\in \mathbb N^+$\\
$\mathrm{(ii)}\>\>\> (x^2.y).x=x^2.(y.x),\>\> \text{for any}\>\>  x,y\in J$
 \end{theorem}
A finite-dimensional real commutative algebra satisfying (\ref{formre}) and the equivalent conditions (i) and (ii) of the theorem is referred to as a finite-dimensional {\sl Euclidean} (or {\sl formally real}) {\sl Jordan algebra}, \cite{ber:2000},  \cite{mcc:2004}, \cite{ior:2009}.\\

Condition (i) of Theorem \ref{POJO} is referred to as {\sl power associativity} while Condition (ii) is referred to as the {\sl Jordan identity}. It is well known and easy to show that the Jordan identity implies the power associativity but it is a very nontrivial result that in the above context the power associativity combined with the formal reality condition imply the Jordan identity \cite{jor-vn-wig:1934}. One can find in \cite{alb:1950} other natural conditions on commutative power-associative algebras leading to the Jordan identity.\\

It is worth noticing that Condition (3.2) and Condition (i) of Theorem~1 (i.e. power associativity) are exactly the conditions needed to have spectral resolutions of the elements of $J$. Indeed a finite-dimensional Euclidean Jordan algebra $J$ has a unit $\bbbone\in J$ such that $\bbbone . x=x$ for any $x\in J$ and any $x\in J$ has a {\sl finite spectral resolution}
\begin{equation}
x=\sum_{r\in I_x}\lambda_r e_r
\label{sprx}
\end{equation}
where $I_x$ is a finite set, $\lambda_r\in \mathbb R$ and $e^2_r=e_r\in J$ for any $r\in I_x$, $e_r . e_s=0$ for $r,s\in I_x$ with $r\not=s$ and 
\begin{equation}
\bbbone=\sum_{r\in I_x}e_r
\label{spr1x}
\end{equation}
i.e. $(e_r)$ is a finite family of orthogonal idempotents in $J$ the sum of which is the unit $\bbbone\in J$. Furthermore the number of elements of $I_x$, $\text{card}(I_x)$, is bounded by a finite number which only depends on $J$ that is one has $\sup_{x\in J}(\text{card}(I_x))<\infty$.\\

This means that one can set
\begin{equation}
x^0=\bbbone
\label{x01}
\end{equation}
and that
\begin{equation}
P\mapsto P(x)=\sum_r P(\lambda_r)e_r
\label{fcalc}
\end{equation}
defines a homomorphism of unital algebras from the algebra $\mathbb R[X]$ of real polynomials into $J$ for any $x\in J$.\\

These finite-dimensional Euclidean Jordan algebras were classified in \cite{jor-vn-wig:1934} where the following theorem is proved.

\begin{theorem}\label{ClassF}
Any finite-dimensional Euclidean Jordan algebra is a direct sum of a finite number of simple ideals. A finite-dimensional simple Euclidean Jordan algebra is isomorphic to one of 
\[
\begin{array}{l}
\mathbb R\ \mathrm{and}\ \ JSpin_{n+2} =\mathbb R\bbbone\oplus \mathbb R^{n+2},\\
\\
J^1_{n+3} = \calh_{n+3}(\mathbb R),
J^2_{n+3} = \calh_{n+3}(\mathbb C),\\
\\
J^4_{n+3} = \calh_{n+3}(\mathbb H),\ \mathrm{with}\ n\in \mathbb N\\
\\
\mathrm{and}\ J^8_3 = \calh_3(\mathbb O).
\end{array}
\]
\end{theorem}

In the above statement $JSpin_{n+1}=\mathbb R\bbbone\oplus \mathbb R^{n+1}$ is endowed with the product
\[
(r\bbbone\oplus v). (r'\bbbone\oplus v')=(rr'+\langle v,v'\rangle)\bbbone\oplus (rv'+r'v)
\]
where $\langle,\rangle$ is the scalar product of the Euclidean space $\mathbb R^{n+1}$ and, for any involutive algebra $\cala$, $\calh_n(\cala)$ denotes the space of hermitian $n\times n$-matrices with entries in $\cala$ endowed with the product
\begin{equation}
A. B = \frac{1}{2} (AB+BA)
\label{anticom}
\end{equation}
where $(A,B)\mapsto AB$ is the product in $M_n(\cala)$.\\

The absence of $JSpin_{1}$ in the above list comes from the fact that $JSpin_{1}$ is not simple but is isomorphic to the commutative associative algebra $\mathbb R\oplus \mathbb R$ of the real functions on a two-points set. Furthermore, notice that one has the following coincidence
\[
\begin{array}{l}
\calh_1(\mathbb R)=\calh_1(\mathbb C)=\calh_1(\mathbb H)=\calh_1(\mathbb O)=\mathbb R,\\
\calh_2(\mathbb R)=JSpin_2,\ \calh_2(\mathbb C)=JSpin_3, \calh_2(\mathbb H)=JSpin_5,\\
\calh_2(\mathbb O)=JSpin_9
\end{array}
\]
and that the $\calh_n(\mathbb O)$ are not Jordan algebras for $n\geq 4$.\\

In each finite-dimensional simple Euclidean Jordan algebra $J$, its unit $\bbbone$ admits a decomposition
\[
\bbbone = \sum^{c(J)}_{i=1} e_i
\]
with $e_i. e_j=0$ for $i\not= j$ (orthogonality) where the $e_i$ are primitive idempotents. The number $c(J)$ does only depend on $J$ and is refered to as the {\sl capacity} (or the {\sl degree}) of $J$. One has $c(\mathbb R)=1$, $c(JSpin_{n+2})=2$, $c(J^1_3)=c(J^2_3)=c(J^4_3)=c(J^8_3)=3$ and $c(J^1_n)=c(J^2_n)=c(J^4_n)=n$ for $n\geq 4$.\\

Note the obvious fact that a finite-dimensional Euclidean Jordan algebra $J$ which is associative is the algebra of real functions on a finite space $K$ which can be identified to the set of {\sl characters} of $J$, that is the set of homomorphisms
\[
\chi : J\rightarrow \mathbb R
\]
of unital algebras. More generally, it will be convenient in the following to consider that a finite-dimensional Euclidean Jordan algebra $J$ is the algebra of ``{\sl real function}" on some associated ``{\sl finite quantum space}" (a dual object). This is justified by the above spectral properties as well as by the fact that $J$ is the set of observable of a finite quantum system (see below). Indeed the pure states on $J$ are the primitive idempotents in $J$ and the transition probabilities between two states are the appropriate traces of the corresponding product of idempotents.\\

Recall finally that a {\sl Jordan algebra} over a commutative field $\mathbb K$ is a $\mathbb K$-vector space endowed with a symmetric bilinear product $(x,y)\mapsto x. y$ satisfying Property (ii) of Theorem \ref{POJO}. For infinite-dimensional real Jordan algebras, there are various generalizations of the formal reality condition (\ref{formre}). In this frame there are (bounded) generalizations of Theorem \ref{POJO}, see \cite{ioc-lou:1985} for the $JB$-algebra version and \cite{vn:1936} for a weak version.

\subsection{The exceptional finite quantum space} 

Apart $J^8_3$, all the other finite-dimensional simple Euclidean Jordan algebras can be realized as real vector spaces of self-adjoint operators in finite-dimensional Hilbert spaces, stable by the anticommutator of the associative composition of linear operators. They can therefore be considered as algebras of observabes of some finite quantum systems.These special algebras enter in the frame of noncommutative geometry based on associative $\ast$-algebras and each nonexceptional simple Lie algebra occurs as Lie algebra of derivations of one of these special simple Euclidean Jordan algebras.\\

However, as proved in \cite{alb:1934} the Jordan algebra $J^8_3$ is {\sl exceptional} which means that it cannot be inbedded in an associative algebra as a subspace stable by the anticommutator of the associative product as Jordan product.\\

Let us recall that if $\cala$ is an associative algebra with product $(x,y)\mapsto xy$, then the anticommutator $(x,y)\mapsto x.y=\frac{1}{2}(xy+yx)$ is a Jordan product on $\cala$ and one denotes by $\cala^{Jord}$ the corresponding Jordan algebra. A Jordan algebra which is isomorphic to a Jordan subalgebra of some $\cala^{Jord}$ as above is said to be {\sl special}, otherwise it is said to be {\sl exceptional}.\\

In spite of the fact that $J^8_3$ is exceptional it can be consistently considered as the algebra of observable of a finite quantum system, see e.g. \cite{jor-vn-wig:1934}, \cite{gun-pir-rue:1978}. Thus it corresponds to an {\sl exceptional quantum space} which is in fact unique \cite{zel:1979}, \cite{zel:1983}. Notice that since $\dim_\mathbb R(J^8_3)=27<\infty$, we consider that the ``corresponding quantum space" is ``finite".\\

The group of all automorphisms of $J^8_3$ is the exceptional Lie group $F_4$ with Lie algebra $\Lie(F_4)=\mathfrak{f}_4$ which is the Lie algebra of all derivations of $J^8_3$. Therefore at first sight  $F_4$ plays the role of diffeomorphism group of the exceptional quantum space while the exceptional simple Lie algebra $\mathfrak{f}_4$ plays the role of the Lie algebra of vector fields on this exceptional quantum space.\\
%

\noindent \underbar{Remark}. Notice that these actions are inner. It follows that $F_4$ and its Lie algebra $\frak f_4$ admit an alternative interpretation in terms of gauge group and Lie algebra. This applies as well for their subgroups and Lie subalgebras. Of particular interest in view of the analysis of Section 2 is the subgroup of $F_4$ which preserves the representation of the octonion algebra $\mathbb O$ as $\mathbb C\oplus\mathbb C^3$, (see in Section 4).

\subsection{Exceptional almost classical quantum spacetime}

As in \S 2.4 these algebras $J^8_3$ over each spacetime point are the fiber of a bundle in algebras $J^8_3$ which is associated to the principal (color) $SU(3)$-bundle on which the gluon field is a connection. This connection is compatible with the identification of the octonionic matrix elements of $J^8_3$ as elements
of $\mathbb C\oplus \mathbb C^3$. Let $J^8_3(M)$ be the real Jordan algebra of the sections of this bundle in $J^8_3$ algebras. The matrix elements of $J^8_3(M)$ are then elements of $A(M)$ (see in \S 2.4). We denote by the same symbol $\nabla_X$ the covariant derivative along the vector field $X\in \Der(C^\infty(M))$ corresponding to the $SU(3)$-connection and by $\rho$ the restriction of $\Der(J^8_3(M))$ to $C^\infty(M)=C^\infty(M)\otimes \bbbone\subset J^8_3(M)$. Again $X\mapsto \nabla_X=s_\nabla(X)$ is a splitting of the exact sequence of $C^\infty(M)$-modules
\begin{equation}
0\rightarrow \text{Int} (J^8_3(M))\stackrel{\subset}{\rightarrow} \Der(J^8_3(M))\stackrel{\rho}{\rightarrow}\Der(C^\infty(M))\rightarrow 0
\label{seqJ}
\end{equation}
where now $\text{Int}(J^8_3(M))=C^\infty(M,\fracf_4)$ with the exceptional Lie algebra $\fracf_4$ identified with the Lie algebra $\Der(J^8_3)$ of derivations of the exceptional Jordan algebra $J^8_3$. The derivations of $J^8_3$ being inner, one has $\Der(J^8_3)=\text{Int}(J^8_3)$ so  $\text{Int}(J^8_3(M))=C^\infty(M,\fracf_4)$. Setting for $\delta \in \Der(J^8_3(M))$
\[
\pi_\nabla(\delta)=\delta-\nabla_{\rho(\delta)}
\]
as in \S 2.4, one defines a projection $\pi_\nabla$ of $\Der(J^8_3(M))$ onto $\text{Int}(J^8_3(M))=C^\infty(M,\fracf_4)=\ker(\rho)$.\\

Assuming as before that $M\cong \mathbb R^4$, one makes the identifications
\begin{equation}
J^8_3(M)=C^\infty(M,J^8_3)=C^\infty(M)\otimes J^8_3
\label{IdJ}
\end{equation}
and, by an abuse of language, one defines the exceptional {\sl almost classical quantum spacetime} by saying that its algebra of real ``smooth functions"  is $J^8_3(M)$. With the identification (\ref{IdJ}), it is the product of the classical spacetime $M$ with the exceptional finite quantum space corresponding to $J^8_3$.  More generally, the notion of almost classical quantum spacetime, or simply of {\sl almost classical spacetime}, is obtained by replacing in (\ref{IdJ}) $J^8_3$ by an arbitrary finite-dimensional Euclidean Jordan algebra.\\

The automorphism group of $J^8_3$ is the exceptional group $F_4$ while the automorphism group of $C^\infty(M)$ is the group $\Diff(M)$ of diffeomorphisms of $M$. The automorphism group $\Aut(C^\infty(M, J^8_3))$ of $C^\infty(M,J^8_3)= C^\infty(M)\otimes J^8_3 $ is the semidirect product $C^\infty(M,F_4) \rtimes \Diff(M)$ of the group $C^\infty(M,F_4)$ with $\Diff(M)$. Thus in a sense this automorphism group plays the role of the diffeomorphism group of the exceptional almost classical quantum spacetime.
Correspondingly the Lie algebra $\Der(C^\infty (M,J^8_3))$ of derivations of the Jordan algebra $C^\infty(M,J^8_3)$ is the semidirect sum of the Lie algebra $C^\infty(M,\mathfrak{f}_4)$ and the Lie algebra $\Der(C^\infty(M))$ of vector fields on $M$ and plays in the same above sense the role of the Lie algebra of vector fields on the exceptional almost classical quantum spacetime.\\

However Section 2.4 and the remark of \S 3.3 suggest to give to $F_4$ and its subgroups an alternative interpretation as structure groups for a gauge theory on $J^8_3$-modules.\\

We wish to develop a field theory on this quantum spacetime. For this one has to introduce the appropriate notion of (bi-)module over $C(M,J^8_3)$, the relevant matter's fields being then elements of such bimodules. It is clear that one needs first a description of modules over $J^8_3$. Furthermore, in view of the previous interpretation of the quark-lepton symmetry, one has to describe the subgroup of $F_4$ acting on such modules which is compatible with the representation of $\mathbb O$ as $\mathbb C\oplus \mathbb C^3$.\\

\section{Modules and fundamental fermions}
\subsection{Modules over Jordan algebras}

There are general definitions of bimodules for classes of algebras \cite{sch:1995}, \cite{jac:1968} or for categories of algebras \cite{mdv:2001}. For Jordan algebras this reduces to the following (see also in \cite{kas-ovs-she:2011}). Let $J$ be a Jordan algebra and let $M$ be a vector space equipped with left and right actions of $J$
\[
J\otimes M\rightarrow M,\>\>\> x\otimes\Phi\mapsto x\Phi\in M
\]
\[
M\otimes J\rightarrow M, \>\>\> \Phi\otimes x\mapsto \Phi x\in M
\]
for $x\in J$, $\Phi\in M$. Define a bilinear product on $J\oplus M$ by setting
\[
(x+\Phi)(x'+\Phi')=xx'+(x\Phi'+\Phi x').
\]
Then $M$ is said to be {\sl a Jordan bimodule over} $J$ if $J\oplus M$ endowed with the above product is a Jordan algebra. This is equivalent to the following conditions (i), (ii) and (iii) :
\[
\left\{
\begin{array}{l}
\mathrm {(i)}\>\>\> x\Phi=\Phi x\\
\mathrm {(ii)}\>\>\> x(x^2\Phi)=x^2(x\Phi)\\
\mathrm {(iii)}\>\>\> (x^2y)\Phi -x^2(y\Phi)=2((xy)(x\Phi)-x(y(x\Phi)))
\end{array}
\right.
\]
for $x,y\in J$ and $\Phi\in M$. When $J$ is unital with unit denoted by $\bbbone$, the Jordan bimodule $M$ will be said to be {\sl unital} whenever it satisfies the further condition (iiii)\\
\phantom{condition (iiii)} $\mathrm{(iiii)}\>\>\>\> \bbbone\Phi=\Phi$\\
for any $\Phi\in M$, which implies that $J\oplus M$ is again a unital Jordan algebra with unit $\bbbone$.\\

It follows from Condition (i) that one can restrict attention to the left action (for instance), this is why Jordan bimodules are usually refered to as {\sl Jordan modules}.\\

Let us define the linear mapping $x\mapsto L_x$ of $J$ into the algebra $\call(M)$ of linear endomorphisms of the $J$-module $M$ by setting
\begin{equation}
L_x\Phi=x\Phi
\label{LM}
\end{equation}
for $x\in J$ and $\Phi\in M$. Then (ii) reads
\begin{equation}
[L_x,L_{x^2}]=0
\label{Lxx2}
\end{equation}
for $x\in J$, while (iii) reads
\[
L_{x^2y}-L_{x^2}L_y-2L_{xy}L_x+2L_xL_yL_x=0
\]
for $x,y\in J$ which is equivalent to
\begin{equation}
L_{x^3}-3L_{x^2}L_x+2L^3_x=0
\label{Lx3}
\end{equation}
and
\begin{equation}
[[L_x,L_y],L_z]+L_{[x,z,y]}=0
\label{intder}
\end{equation}
for $x,y,z\in J$ where $[x,z,y]=(xz)y-x(zy)$ is the associator.\\

Thus Conditions (ii) and (iii) can be replaced by Conditions (\ref{Lxx2}), (\ref{Lx3}) and (\ref{intder}) above while the unital condtion (iiii) reads 
\begin{equation}
L_{\bbbone} =\bbbone
\label{L1}
\end{equation}
where on the right-hand side $\bbbone$ denotes the unit of the algebra $\call(M)$.\\

Notice that if $p\in J$ is an idempotent, i.e. satisfies $p^2=p$, then (\ref{Lx3}) reads 
\[
L_p-3 L^2_p+ 2L^3_p=0
\]
that is
\begin{equation}
L_p(L_p-\frac{1}{2})(L_p-1)=0
\label{p}
\end{equation}
from which follows the associated Pierce decomposition of the module $M$.

\subsection{$J^8_3$-modules and $SU(3)\times SU(3)$-action}

Any Jordan algebra $J$ is canonically a Jordan module over itself which is unital whenever $J$ has a unit. The list of the unital irreducible Jordan modules over the finite-dimensional simple Euclidean Jordan algebras can be found in \cite{jac:1968}.
In the case of the exceptional Jordan algebra $J^8_3$, it turns out that any unital Jordan module is a product or a direct sum of modules isomorphic to $J^8_3$, 
\cite{jac:1968}. In particular any finite module $M$ over $J^8_3$ is a finite product (or direct sum) of modules isomorphic to $J^8_3$ so that one has $M=J^8_3\otimes E$ for some finite-dimensional real vector space $E$.\\

In view of our interpretation of the quark-lepton symmetry it is important to describe the subgroup of $F_4$ which preserves the representations of the octonions occuring in the elements of $J^8_3$ as elements of $\mathbb C\oplus\mathbb C^3$. For this, following  \cite{yok:2009}, one associates to the element
\begin{equation}
\left(
\begin{matrix}
\zeta_1 & x_3 & \bar x_2\\
\bar x_3 & \zeta_2 & x_1\\
x_2 & \bar x_1 & \zeta_3
\end{matrix}
\right)
\label{jel}
\end{equation}
of $J^8_3$ the following element
\begin{equation}
\left(
\begin{matrix}
\zeta_1 & z_3 & \bar z_2\\
\bar z_3 & \zeta_2 & z_1\\
z_2 & \bar z_1 & \zeta_3
\end{matrix}
\right) +(Z_1,Z_2,Z_3)
\label{cjl}
\end{equation}
of $J^2_3\oplus M_3(\mathbb C)$ where
\begin{equation}
x_i=z_i+Z_i\in \mathbb C \oplus \mathbb C^3
\label{coc}
\end{equation}
are the representation in $\mathbb C\oplus \mathbb C^3$ of the 3 elements $x_i$ of $\mathbb O$. The desired subgroup of $F_4$ is the group $SU(3)\times SU(3)/\mathbb Z_3$ with action induced by the action of $SU(3)\times SU(3)$ on $J^2_3\oplus M_3(\mathbb C)$ given by \cite{yok:2009}
\begin{equation}
H\mapsto VHV^\ast, \>\>\> M\mapsto UMV^\ast
\label{acg}
\end{equation}
for $(U,V)\in SU(3)\times SU(3)$ and $H\in J^2_3=\calh_3(\mathbb C), M\in M_3(\mathbb C)$. Notice that the action of the first factor $U$ is the one which corresponds to the previous action of $SU(3)$ on $\mathbb C\oplus \mathbb C^3$ $(\simeq \mathbb O)$ described in Section 2. Thus with the previous interpretation, the first factor $SU(3)$ is the ``colour group". The second factor $SU(3)$ mixes the generations and furthermore mixes the complexes $z_i$ corresponding to generations of leptons with the reals $\zeta_i$ of the diagonal which correspond in this picture to new particles of spin 1/2. Since the mass scales of the 3 generations are very different this latter factor $SU(3)$ must be strongly broken, however there could remain some unbroken finite subgroup of this $SU(3)$ which should then play a role in this approach for the see-saw mechanism and the structure of the fermion mass matrices.

 \subsection{Modules for the two families by generation}
 
 As recalled in \S 3.1, there are two families by generation. In view of \S 4.2 it seems, at first sight, reasonable to take as $J^8_{3}$-module the product of two copies $J^u$ and $J^d$ of $J^8_3$ as module with the particle assignment
 \begin{equation}
 J^u=\left ( \begin{array}{ccc}
 \alpha_1 & \nu_\tau+t & \bar\nu_\mu-c\\
 \bar\nu_\tau-t & \alpha_2 & \nu_e+u\\
 \nu_\mu+c & \bar\nu_e-u  &\alpha_3
 \end{array}
 \right)
 \label{Ju1}
 \end{equation}
 
 \begin{equation}
 J^d=\left ( \begin{array}{ccc}
 \beta_1 & \tau+b & \bar\mu-s\\
 \bar\tau-b & \beta_2 & e+d\\
 \mu+s & \bar e-d & \beta_3
 \end{array}
 \right)
 \label{Jd1}
 \end{equation}
 
or the representation (\ref{cjl})
 
 \begin{equation}
 J^u=\left( \begin{array}{ccc}
 \alpha_1 & \nu_\tau & \bar\nu_\mu\\
 \bar\nu_\tau & \alpha_2 & \nu_e\\
 \nu_\mu & \bar\nu_e & \alpha_3
 \end{array}
 \right)
 + (u,c,t)
 \label{Ju2}
 \end{equation}
 
 \begin{equation}
 J^d=\left( \begin{array}{ccc}
 \beta_1 & \tau & \bar\mu\\
 \bar\tau & \beta_2 & e\\
 \mu & \bar e & \beta_3
 \end{array}
 \right)
 + (d,s,b)
 \label{Jd2}
 \end{equation}
 where we have identified the fundamental fermions with there internal spaces ($\mathbb C^3$ for the quarks and $\mathbb C$ for the leptons). The diagonals correspond then to new spin1/2 fermions with $\mathbb R$ as internal spaces. It is therefore natural to consider these particles $\alpha_k, \beta_k$ ($k\in\{1,2,3\}$) as described by Majorana spinors. The interest of this latter identifications is that the addition of such particles does not spoil the anomalies cancellation of the Standard Model. In this context the corresponding module over $C^\infty(M,J^8_3)$ would be the tensor product of the $C^\infty(M)$-module of Majorana spinor fields over spacetime by the above $J^8_3$-module.\\
 
 \subsection{Problem of the $U(1)\times SU(2)$-symmetry}
 
 Up to now the analysis deals only with the color symmetry $SU(3)$, its tentative connection with the quark-lepton symmetry and the 3 generations. This leads naturally to the above assignment of internal spaces for fundamental fermions (of \S 4.3) by (\ref{Ju2}) and (\ref{Jd2}) together with an action of $SU(3)\times SU(3)/\mathbb Z_3$ described in \S 4.2 by (\ref{acg}) where the first factor $SU(3)$ is the color group.\\
 
 However then, the description of the $U(1)\times SU(2)$-symmetry may become problematic. Either one can solve the problem by $\gamma_5$-combinatorics or one should work with a product of four factors $J^u_L, J^u_R, J^d_L, J^d_R$. But in this latter option the cancellation of anomalies could become a non-obvious task.

 In view of the above framework, there is another natural way to explore which consists in adding over each spacetime point to the exceptional finite quantum space another finite quantum space in which the group $U(1)\times SU(2)$ is involved. That is to add to $J^8_3$ another finite-dimensional Euclidean Jordan algebra which admits $U(1)\times SU(2)$ as automorphism group preserving some additional structure (like $SU(3)\times SU(3)/\mathbb Z_3$ for $J^8_3$). There is such an algebra, namely the simple Euclidean Jordan algebra
 \[
 J^4_2=\calh_2(\mathbb H)=J \text{Spin} (5)
 \]
 of the hermitian $2\times 2$-matrices with entries in the field of quaternions $\mathbb H$. Let us explain this.\\
 
 A quaternion $q$ can be viewed as two complexes $z_1$ and $z_2$
 \begin{equation}
 q=(z_1,z_2)=z_1+z_2 j
 \label{q2c}
 \end{equation}
 and the subgroup of $\Aut(\mathbb H)=SU(2)/\mathbb Z_2$ which is complex linear, i.e. which preserves the imaginary $i$, is the group $U(1)$ acting as $(z_1,z_2)\mapsto (z_1,e^{i\theta} z_2)$, that is
 \begin{equation}
 z_1+z_2 j\mapsto z_1+e^{i\theta} z_2 j
 \label{U1H}
 \end{equation}
 for $e^{i\theta}\in U(1)$. Thus $U(1)$ is the analog for $\mathbb H$ to $SU(3)$ for $\mathbb O$.\\
 
 Let now consider an element of $J^4_2$ written in terms of $q=z_1+z_2 j$
 \[
 \left(
 \begin{array}{cc}
 \xi_1 & q\\
 \bar q & \xi_2
 \end{array}
 \right) =
 \left(
 \begin{array}{cc}
 \xi_1 & z_1+z_2j\\
 \bar z_1-z_2j & \xi_2
 \end{array}
 \right) =
 \left(
 \begin{array}{cc}
 \xi_1 & z_1\\
 \bar z_1 & \xi_2
 \end{array}
 \right )
 + z_2 \left(
 \begin{array}{cc}
 0 & j\\
 -j & 0
 \end{array}
 \right)
 \]
 where $\xi_1,\xi_2\in \mathbb R$. Let $U\in SU(2)$ be an arbitrary unitary $2\times 2$ complex matrix of determinant 1. By using $jz=\bar z j$ for $z\in \mathbb C$ one gets
 \begin{equation}
 U\left(
 \begin{array}{cc}
 0 & j\\
 -j & 0
 \end{array}
 \right)
 U^\ast=\left(
 \begin{array}{cc}
 0 & j\\
 -j & 0\end{array}
 \right)
 \label{jinv}
 \end{equation}
 from which it follows that the action of $(e^{i\theta},U) \in U(1)\times SU(2)$
 \[
 \left(
 \begin{array}{cc}
 \xi_1 & z_1\\
 \bar z_1 & \xi_2
 \end{array}
 \right) + z_2\left(
 \begin{array}{cc}
 0 & j\\
 -j & 0
 \end{array}
 \right)
 \mapsto U \left(
 \begin{array}{cc}
 \xi & z_1\\
 \bar z_1 & \xi_2
 \end{array}
 \right)
U^\ast +
 e^{i\theta}z_2 \left(
 \begin{array}{cc}
 0 & j\\
 -j & 0
 \end{array}
 \right)
 \]
 is an automorphism of $J^4_2$ which preserves the representation of the quaternions occurring in the elements of $J^4_2$ as elements of $\mathbb C\oplus \mathbb C$. Furthermore any such an automorphism is of this form. Thus $U(1)$ is the analog for $J^4_2$ of the color $SU(3)$ for $J^8_3$ while $SU(2)$ is the analog for $J^4_2$ of the other $SU(3)$ for $J^8_3$ (with the action described by (4.4) in \S 4.2). Notice that the $SU(2)$-action is an action of $SU(2)/\mathbb Z_2=SO(3)$. \\
 
 It is too early to know whether the addition of the factor $J^4_2$ is relevant but if it occurs it is tempting to add the factor $J^2_1=\mathbb R$ since then we would have
 \begin{equation}
 J=\oplus^3_{k=1} J^{2^k}_k
 \label{Jtent}
 \end{equation}
 as finite-dimensional Euclidean Jordan algebra which looks like the semi-simple part of some more elaborated object. Notice that $J$ is the real Jordan algebra of hermitian elements of 
 \begin{equation}
 \mathbb C\oplus M_2(\mathbb H)\oplus M_3(\mathbb O)
 \label{MJ}
 \end{equation}
 which is not associative in view of the occurrence of the last term $M_3(\mathbb O)$.

  \subsection{Charge conjugation}
 In the identification of the octonion algebra $\mathbb O$ as $\mathbb C\oplus \mathbb C^3$, the octonionic conjugate of $(z,Z)=z+Z$ is $\overline{(z,Z)}=(\bar z,-Z)=\bar z-Z$.\\
 On the other hand if one interprets $Z\in \mathbb C^3$ as an internal quark state and $z\in \mathbb C$ as an internal state of the corresponding lepton, the charge conjugation $\calc$ induces the complex conjugation in $\mathbb C^3$ and in $\mathbb C$ that is 
 \begin{equation}
 \calc(z+Z)=\bar z+\bar Z
 \end{equation}
 for $z\in \mathbb C$, $Z\in \mathbb C^3$.\\
 
 It is easy to verify by using Formula (\ref{proc}) that $\calc$ is an automorphism of $\mathbb O$  $(\calc\in G_2)$ with $\calc^2=\bbbone$ by construction. This automorphism of $\mathbb O$ induces canonically an involutive automorphism of the Jordan algebra $J^8_3$ as well as of the $J^8_3$-modules. These automorphisms will be all denoted by $\calc$.\\

 \section{Differential calculus over Jordan algebras}
 
 \subsection{First order differential calculi}
 
 Given a category $\mathbf C$ of algebras over some field $\mathbb K$ and an algebra $\cala$ in $\mathbf C$, a pair $(\Omega^1,d)$ of an $\cala$-bimodule $\Omega^1$ for $\mathbf C$ \cite{mdv:2001}
  and a derivation $d$ of $\cala$ into $\Omega^1$, that is a $\mathbb K$-linear mapping satisfying (the Leibniz rule)
  \[
  d(xy)=d(x)y+xd(y)
  \]
  for any $x,y\in \cala$, will be refered to as a {\sl first order differential calculus over $\cala$} for $\mathbf C$ or simply a first order calculus over $\cala$ when no confusion arises. There is an obvious notion of homomorphism of first order differential calculi over $\cala$. This terminology as well as the one for higher order differential calculi was introduced in \cite{wor:1989} for the case of the category of associative algebras on a field $\mathbb K$.
  
  In the following we shall be concerned with the category of real unital Jordan algebras so, if  $J$ is such a unital Jordan algebra, a $J$-bimodule is a unital $J$-module which will be simply refered to as a $J$-module when no confusion arises.\\
  
  Let $J$ be a unital Jordan algebra and let us define {\sl its center} $Z(J)$ by
  \begin{equation}
  Z(J)=\{z\in J \vert [x,y,z]=[x,z,y]=0, \forall x,y\in J\}
  \label{center}
  \end{equation}
  where $[x,y,z]=(xy)z-x(yz)$ denotes the associator of $x,y,z\in J$. The center $Z(J)$ is a unital associative subalgebra of $J$. The Lie algebra $\Der(J)$ of all derivations of $J$ into itself is also a $Z(J)$-module and we denote by $\Omega^1_{\der}(J)$ the $J$-module of all the $Z(J)$-homomorphisms of $\Der(J)$ into $J$. One defines a derivation
  \[
  d_{\der}:J\rightarrow \Omega^1_{\der}(J)
  \]
  by setting
  \[
  d_{\der}(x)(X)=X(x)
  \label{dDer}
  \]
  for any $x\in J$ and $X\in\Der(J)$.\\
   This first-order differential calculus $(\Omega^1_{\der}(J),d_{\der})$ will be referred to as the {\sl derivation-based first order differential calculus} over $J$.
  
  \subsection{First order differential calculi over $J^8_3$}
  
  The Lie algebra $\Der(J^8_3)$  of all derivations of the exceptional Jordan algebra $J^8_3$ is the exceptional compact simple Lie algebra $\mathfrak f_4$. Since $J^8_3$ is simple $Z(J^8_3)=\mathbb R$ so the unital $J^8_3$-module $\Omega^1_{\der}(J^8_3)$ is the vector space of all linear mappings of $\Der(J^8_3)$ into $J^8_3$ that is 
  \[
  \Omega^1_{\der}(J^8_3)=J^8_3\otimes \Der(J^8_3)^\ast
  \label{omegad}
  \]
  where $\Der(J^8_3)^\ast=\mathfrak f_4^\ast$ denotes the dual vector space of $\Der(J^8_3)$.\\
  
  One defines as before the derivation $d_{\der}:J^8_3\rightarrow \Omega_{\der}(J^8_3)$ by setting $d_{\der}(x)(X)=X(x)$ for any $x\in J^8_3$ and $X\in \Der(J^8_3)$. The first order differential calculus $(\Omega_{\der}(J^8_3),d_{\der})$  over $J^8_3$ is characterized by the following universal property.
  
 \begin{proposition}\label{univ}
 For any first order calculus $(\Omega^1,d)$ over $J^8_3$, there is a unique (unital) Jordan module homomorphism $i_d:\Omega^1_{\der}(J^8_3)\rightarrow \Omega^1$ such that one has $d=i_d\circ d_{\der}$.
 \end{proposition}
 
 \noindent \underline{Proof}. From the general structure of the $J^8_3$-modules, $\Omega^1$ is of the form $\Omega^1=J^8_3\otimes E$ for some vector space $E$. Let $(e^\alpha)$ be a basis of $E$ then $d(x)=X_\alpha(x)\otimes e^\alpha$ with $X_\alpha\in \Der(J^8_3)$ for any $\alpha$. Let $(\partial_k)$ be a basis of $\Der(J^8_3)$ with dual basis $\theta^k$ then one has $d(x)=\partial_k(x)\otimes C^k_\alpha e^\alpha$ with the $C^k_\alpha\in \mathbb R$. Define $i_d$ by $i_d(x\otimes \theta^k)=x\otimes C^k_\alpha e^\alpha$, then $i_d$ is a $J^8_3$-module homomorphism satisfying $d=i_d\circ d_{\der}$ which is clearly unique.~$\square$\\
 
 \noindent \underline{Remark}. Proposition \ref{univ} is very specific to the exceptional Jordan algebra $J^8_3$ and is a direct consequence of the fact that $J^8_3$ is, as $J^8_3$-module, the only irreducible module. For instance if $J$ is a finite-dimensional simple Euclidean Jordan algebras of Theorem 2 which is distinct of $\mathbb R$ and of $J^8_3$, the first order differential calculus $(\Omega^1_{\der}(J), d_{\der})$ is not universal since then there several inequivalent irreducible $J$-modules and not only $J$ itself. In all these cases one has $\Omega^1_{\der}(J)=J\otimes \Der(J)^\ast$ and  $\Der(J)$ is a real compact form of a classical simple Lie algebra. Furthermore all the classical simple Lie algebras are realized in this way as Lie algebras of derivations of the finite-dimensional simple special Euclidean Jordan algebras.
  
 \subsection{Differential graded Jordan algebras}
 
Let us define a {\sl differential graded Jordan algebra} to be a $\mathbb N$-graded algebra
\[
\Omega=\oplus_{n\in \mathbb N}\Omega^n
\]
which is a {\sl Jordan superalgebra}, (see e.g. in \cite{kac:1977}, \cite{mar-zel:2010}), for the induced $\mathbb Z/2\mathbb Z$ degree and which is equipped with a {\sl differential}, that is with an antiderivation $d$ of degree 1 and of square 0.\\

Thus $\Omega$ is a graded-commutative algebra, that is one has
\begin{equation}
ab=(-1)^{\vert a\parallel b\vert} ba \in \Omega^{\vert a\vert+\vert b\vert}
\label{grcom}
\end{equation}
for $a\in \Omega^{\vert a\vert}$, $b\in \Omega^{\vert b\vert}$, and one has the {\sl graded Jordan identity}
\begin{equation}
(-1)^{\vert a\parallel c\vert}[L_{ab},L_c]_{\gr}+(-1)^{\vert b\parallel a\vert}[L_{bc},L_a]_{\gr}+(-1)^{\vert c\parallel b\vert}[L_{ca},L_b]_\gr=0
\label{grJ}
\end{equation}
for $a\in \Omega^{\vert a\vert},b\in\Omega^{\vert b\vert}, c\in \Omega^{\vert c\vert}$ where $L_a$ is the left-multiplication operator by $a\in \Omega$ defined by $L_a(x)=ax$ for any $x\in \Omega$ and where $[\bullet, \bullet]_{\gr}$ denotes the graded commutator
\[
[A,B]_{\gr}=AB-(-1)^{\vert a \parallel b\vert}BA
\]
for $A$ of degree $\vert a\vert$ and $B$ of degree $\vert b\vert$. The differential $d$ satisfies
\[
d^2=0
\]
\[
d\Omega^n\subset \Omega^{n+1}
\]
for any $n\in \mathbb N$, and the graded Leibniz rule 
\[
d(ab)=d(a)b+(-1)^{\vert a\vert} ad(b)
\]
for $a\in \Omega^{\vert a\vert}$, $b\in \Omega$.\\

In the following the real unital differential graded Jordan algebras will be our models for the ``algebras of differential forms" on ``quantum spaces".\\

Let $J$ be a Jordan algebra. A differential graded Jordan algebra $\Omega=\oplus_{n\geq 0}\Omega^n$ such that $\Omega^0=J$ will be called a {\sl differential calculus over} $J$.\\

\subsection{Higher order derivation-based differential calculi}

Let $J$ be a unital Jordan algebra with center $Z(J)$ and let $\Omega^n_{\der}(J)$ be the $J$-module of all $n$-$Z(J)$-linear antisymmetric mappings of $\Der(J)$ into $J$, that is $\omega\in \Omega^n_{\der}(J)$ is a $Z(J)$-linear mapping
\[
\omega:\wedge^n_{Z(J)} \Der(J)\rightarrow J
\]
of the n-th exterior power over $Z(J)$ of the $Z(J)$-module $\Der(J)$ into $J$ as a $Z(J)$-module. Then $\Omega_{\der}(J)=\oplus_{n\geq 0} \ \Omega^n_{\der}(J)$ is canonically a differential graded Jordan algebra with differential given by the Chevalley-Eilenberg formula
\begin{equation}
\begin{array}{ll}
d\omega(X_0,\cdots,X_n) & =  \sum_{0\leq k\leq n} (-1)^k X_k\  \omega(X_0,\stackrel{k\atop \vee}{\cdots},X_n) \\
&+ \sum_{0\leq r<s\leq n}(-1)^{r+s}\ \omega([X_r,X_s],X_0,\stackrel{r\atop \vee}{\cdots} \stackrel{s\atop \vee}{\cdots},X_n)
\end{array}
\label{C-E}
\end{equation}
 for $\omega\in \Omega^n_{\der}(J)$ and $X_p\in \Der(J)$. Thus $\Omega_{\der}(J)$ is a differential calculus over $J$ which will be referred as the {\sl derivation-based differential calculus} over $J$.\\
 
 If $J$ is a finite-dimensional simple Euclidean Jordan algebra, (i.e. one of the list of Theorem 2), one has
 
\begin{equation}
\Omega_{\der}(J)=J\otimes \wedge \Der(J)^\ast
\label{omegaJ}
\end{equation}
where $\wedge\Der(J)^\ast$ is the exterior algebra of the dual of the finite-dimensional real Lie algebra $\Der(J)$.  In the case where $J$ is the exceptional Jordan algebra $J^8_3$, this differential calculus is characterized by the following universal property.
 
 \begin{proposition}\label{dex}
 Any homomorphism $\varphi$ of unital Jordan algebra of $J^8_3$ into the Jordan subalgebra $\Omega^0$ of a unital differential graded Jordan algebra $\Omega=\oplus \Omega^n$ has a unique extension $\tilde\varphi:\Omega_{\der}(J^8_3)\rightarrow \Omega$ as a homomorphism of differential graded Jordan algebras.
 \end{proposition}
 
 \noindent
 \underbar{Proof} (sketch of). By using the fact that the $J^8_3$-modules are of the form $J^8_3\otimes \cale$ for some vector spaces $\cale$ and the definition of the $\mathbb N$-graded Jordan superalgebras, Formulae (\ref{grcom}) and (\ref{grJ}), it follows that any unital differential graded Jordan algebra $\Omega=\oplus \Omega^n$ which contains $J^8_3$ as unital Jordan subalgebra of $\Omega^0$ is as algebra of the form
 \[
 \Omega=J^8_3\otimes \Omega_{(0)}
 \]
 where $\Omega_{(0)}=\oplus \Omega^n_{(0)}$ is an associative graded-commutative algebra (this generalizes a classical result in the non graded case). The proof follows then from Proposition \ref{univ},  from $d^2=0$ with the graded Leibniz rule for $d$ and from the simplicity  of $J^8_3$.~$\square$\\
 
 A more detailed proof will appear in a forthcoming paper.\\
 
 \noindent
 \underbar{Remark}. Notice that as Proposition \ref{univ},  Proposition \ref{dex} is very specific to the exceptional Jordan algebra $J^8_3$.\\
 
 As a consequence all differential calculus over $J^8_3$ which are generated as differential graded Jordan algebra by $J^8_3$ are quotients of $\Omega_\der(J^8_3)$, in particular among these, one has all the differential graded Jordan algebras of the form
 
 \begin{equation}
 \Omega=J^8_3\otimes \wedge \fracg^\ast
 \label{diffexc}
 \end{equation}
 where $\fracg$ is a Lie subalgebra of $\Der(J^8_3)=\frak f_4$, the differential being the Chevalley-Eilenberg differential ($\fracg$ acting by derivations on $J^8_3$).
 
 \subsection{Differential calculi as $A_\infty$-algebras}
 
 Let $\cala$  be a complex unital $\ast$-algebra and  let $\Omega_u(\cala)=\oplus_{n\geq 0}\Omega^n_u(\cala)$ be the universal differential calculus over $\cala$, \cite{kar:1983}. It is shown in \cite{wor:1989} (see also in \cite{mdv:2001}) that there is a unique involution on $\Omega_u(\cala)$ which extends the involution of $\cala$ for which it is a differential graded $\ast$-algebra. This means that one has
 \begin{equation}
 \left\{
 \begin{array}{l}
 (\alpha\beta)^\ast=(-1)^{ab}\beta^\ast \alpha^\ast\\
 d(\gamma^\ast)=(d\gamma)^\ast
 \end{array}
 \right.
 \label{stardga}
 \end{equation}
 for any $\alpha\in\Omega^a_u(\cala),\beta\in\Omega^b_u(\cala)$ and $\gamma\in\Omega_u(\cala)$. The cohomology of $\Omega_u(\cala)$ is trivial, i.e. $H^n(\Omega_u(\cala))=0$ for $n\geq 1$ and $H^0(\Omega_u(\cala))=\mathbb C$, and in fact by using a linear form $\omega$ on $\cala$ such that $\omega(\bbbone)=1$, one constructs a linear mapping $K$ of degree -1 on $\Omega_u(\cala)$ such that
 \begin{equation}
 dK+Kd=I
 \label{coho}
 \end{equation}
 on $\oplus_{n\geq 1}\Omega^n_u(\cala)$ (see in \cite{mdv:2001}). By replacing $K$ by $\frac{1}{2}(K+K^\ast)$ one can assume that $K(\alpha^\ast)=K(\alpha)^\ast$ for any $\alpha\in\Omega_u(\cala)$.\\
 
 Let $J=\calh(\cala)$ be the real Jordan algebra of the hermitian elements of $\cala$, i.e. $J=\{h\in \cala\vert h^\ast=h\}$ endowed with the product
 \[
 h\circ h'=\frac{1}{2}(hh'+h'h)
 \]
 for $h,h'\in J$. Then the real graded subspace $\Omega(J)$ of $\Omega_u(\cala)$ of hermitian elements is a real differential graded Jordan algebra, that is a differential calculus over $J$ for the graded Jordan product defined by
 \begin{equation}
 \alpha\circ \beta=\frac{1}{2}(\alpha\beta+(-1)^{ab}\beta\alpha)
 \label{grJ}
 \end{equation}
 for $\alpha\in\Omega^a(J),\beta\in\Omega^b(J)$.\\
 
 The contracting homotopy $K$ restricts to $\Omega(J)$. By using $K$, one constructs by induction on $n$ a sequence of product $m_n(\alpha_1,\cdots,\alpha_n)$ of respective degrees $2-n$ on $\Omega(J)$ starting with $m_1(\alpha)=d\alpha$ and $m_2(\alpha,\beta)=\alpha\circ \beta$ such that, endowed with these products, $\Omega(J)$ is an $A_\infty$-algebra \cite{kel:2001}. This implies in particular that the graded Jordan product $\alpha\circ \beta$ of $\Omega(J)$ is associative up to homotopy. In fact, by taking
 \begin{equation}
 m_3(\alpha,\beta,\gamma)=K((\alpha\circ\beta)\circ \gamma-\alpha\circ(\beta\circ \gamma))
 \label{m3}
 \end{equation}
 for $\alpha\in\Omega^a(J), \beta\in\Omega^b(J),\gamma\in\Omega^c(J)$ with $a+b+c\geq 1$, one gets 
 \[
 (\alpha\circ\beta)\circ \gamma-\alpha\circ (\beta\circ \gamma)=d(m_3) (\alpha,\beta,\gamma)
 \]
 i.e. the associativity up to homotopy of the graded Jordan product.\\
 
 The occurrence of $A_\infty$-structures  in the present context is not accidental. For instance it will be shown in another paper that, for the simple finite dimensional Euclidean Jordan algebras of Theorem 2, the derivation-based differential calculus admits such $A_\infty$-structures with $m_1$ given by the differential and $m_2$ given by the graded Jordan product.

 \section{Connections on Jordan modules}
 
 \subsection{Derivation-based connections. First definition}
 
 Let $J$ be a unital Jordan algebra with center (= centroid) $Z(J)$ and let $\Der(J)$ be the Lie algebra and $Z(J)$-module of derivations of $J$. A {\sl (derivation-based) connection} $\nabla$ on a unital $J$-module $M$ is a linear mapping $X\mapsto \nabla_X$ of $\Der(J)$ into the linear endomorphisms of $M$ such that
 \begin{equation}
 \left\{
 \begin{array}{l}
 \nabla_X(xm)=X(x)m +x\nabla_X(m)\\
 \\
 \nabla_{zX}(m)=z\nabla_X(m)
 \end{array}
 \right.
  \label{covder}
 \end{equation}
 for any $m\in M$, $x\in J$ and $z\in Z(J)$.\\
 
 It follows from this definition that the difference $\nabla-\nabla'$ between 2 connections on $M$ is a $Z(J)$-linear mapping of $\Der(J)$ into the $Z(J)$-module of all the $J$-module endomorphisms of $M$.\\
 
It also follows that 
\begin{equation}
R_{X,Y}=[\nabla_X,\nabla_Y]-\nabla_{[X,Y]}
\end{equation}
 satisfies
 \begin{equation}
 \left\{
 \begin{array}{l}
 R_{X,Y}(xm)=x R_{X,Y}(m)\\
 \\
 R_{zX,Y}(m)=z R_{X,Y}(m)
 \end{array}
 \right.
 \label{curv}
 \end{equation}
 in other words that $R$ is a $Z(J)$-linear mapping of $\wedge^2_{Z(J)}\Der(J)$ into the $Z(J)$-module of all the $J$-module endomorphism of $M$ which will be referred to as the {\sl curvature} of $\nabla$.\\
 
 It is clear that if $\fracg$ is a Lie subalgebra and a $Z(J)$-submodule of $\Der(J)$, the restriction of (\ref{covder}) to $\fracg$, i.e. for $X\in\fracg$, makes sense and the corresponding notion will be referred to as a {\sl derivation-based $\fracg$-connection} on $M$.
 
 \subsection{Derivation-based connections. Second definition}
 
 Let $J$, $Z(J)$, $\Der(J)$ and $M$ be as in \S 6.1 and let $\Omega^n_{\der}(M)$ be the $J$-module of all $n$-$Z(J)$-linear antisymmetric mapping of $\Der(J)$ into $M$, which means that  $\Phi\in \Omega^n_{\der}(M)$ is a $Z(J)$-linear mapping
 \[
 \Phi: \wedge^n_{Z(J)}\Der(J)\rightarrow M
 \]
 with the notations of \S 5.4. The graded $J$-module $\Omega_\der(M)=\oplus_{n\geq 0}\Omega^n_\der(M)$ in naturally a graded Jordan module over the graded Jordan algebra $\Omega_\der(J)$: The product of $\omega\in \Omega^m_\der(J)$ with $\Phi\in\Omega^n_\der(M)$ is the element $\omega\Phi\in\Omega^{m+n}_\der(M)$ obtained by product of evaluations on derivations followed by antisymmetrization in the derivations. A {\sl (derivation-based) connection} on $M$ is a linear endomorphism $\nabla$ of $\Omega_\der(M)$ such that
 \begin{equation}
 \left\{
 \begin{array}{l}
 \nabla(\Omega^n_\der(M))\subset \Omega^{n+1}_\der(M)\\
 \\
 \nabla(\omega \Phi)=d(\omega) \Phi + (-1)^m\omega\nabla (\Phi)
 \label{covdiff}
 \end{array}\right.
 \end{equation}
 for any $m, n\in \mathbb N$, $\omega\in \Omega^m_\der(J)$ and $\Phi\in \Omega_\der(M)$. Let $\nabla$ be such a connection and define $\nabla_Xm$ by
 \begin{equation}
 \nabla(m)(X)
 \label{cocod}
 \end{equation}
 for $m\in M=\Omega^0_\der(M)$ and $X\in \Der(J)$. Then $X\mapsto \nabla_X$ is a connection on $M$ in the sense of \S 6.1. Conversely if $X\mapsto \nabla_X$ is a connection in the sense of \S 6.1 one defines $\nabla$ on $\Omega^n_\der(M)$ by
 \begin{equation}
 \begin{array}{ll}
 \nabla(\Phi)(X_0,\cdots,X_n) & = \sum^n_{p=0}(-1)^p\nabla_{X_p}(\Phi(X_0,\stackrel{p\atop \vee}{\cdots\cdots}, X_n))\\
 &+ \sum_{0\leq r<s\leq n}(-1)^{r+s}\ \Phi([X_r,X_s],X_0,\stackrel{r\atop \vee}{\cdots} \stackrel{s\atop \vee}{\cdots},X_n)
 \end{array}
 \label{evcov}
 \end{equation}
 for $\Phi\in \Omega^n_\der(M)$, $X_p\in\Der(J)$. This $\nabla$ satisfies the axioms (\ref{covdiff}). We shall refer to $\nabla$ as the {\sl covariant differential} while $\nabla_X$ for $X\in \Der(J)$ will be refered to as the covariant {\sl derivative along} $X$.\\
 
 From the axioms (\ref{covdiff}), it follows that one has
 \begin{equation}
 \nabla^2(\omega\Phi)=\omega\nabla^2(\Phi)
 \label{covcurv}
 \end{equation}
 for any $\omega\in \Omega_\der(J)$ and $\Phi\in\Omega_\der(M)$ furthermore
 \[
 \nabla^2(m)(X,Y)=R_{X,Y}(m)
 \]
 for $m\in M$. Thus $\nabla^2$ is an homomorphism of $\Omega_\der(J)$-module which is also called {\sl the curvature} of $\nabla$ since it is expressible in terms of the $R_{X,Y}$ and conversely.
 
 \subsection{General connections}
 
 Let $\Omega=\oplus_{n\geq 0}\Omega^n$ be a differential graded Jordan algebra and let $\Gamma=\oplus_{n\geq 0}\Gamma^n$ be a graded Jordan module over the graded Jordan algebra $\Omega$, (the axioms for this notion are easy to guess). A {\sl connection} on $\Gamma$ will be defined to be a linear endomorphism $\nabla$ of $\Gamma$ satisfying 
 \begin{equation}
 \left\{
\begin{array}{l}
\nabla(\Gamma^n)\subset \Gamma^{n+1}\\
\nabla(\omega\Phi)=d(\omega)\Phi+(-1)^m\omega\nabla(\Phi)
\label{gcovd}
\end{array}
\right.
\end{equation}
for $m, n\in \mathbb N$, $\omega\in\Omega^m$ and $\Phi\in \Gamma$.\\

The axioms (\ref{gcovd}) imply that one has
\begin{equation}
\nabla^2(\omega\Phi)=\omega\nabla^2(\Phi)
\label{gcurv}
\end{equation}
for any $\omega\in\Omega$ and $\Phi\in \Gamma$, so $\nabla^2$ is an homogeneous $\Omega$-module homomorphism of degree 2 which will be refereed to as the {\sl curvature} of the connnection $\nabla$.\\

Notice that the formalism of 6.2 is a particular case of this formalism and that within it, the Bianchi identity reduces to the trivial identity
\[
\nabla\nabla^2=\nabla^2\nabla
\]
that is to the associativity of the composition of the endomorphism $\nabla$.

\section{Tentative conclusion}

From a physical point of view, it is clear that what is described in these notes is quite incomplete : One should write some dynamics. Before that,  one must develop several points.\\

Firstly one must understand the formulation of the $U(1)\times SU(2)$-symmetry in this frame. Since the new suspected particles have $\mathbb R$ as internal spaces, it is suggested that they are described by Majorana fermions in a 4-Lorentzian approach where the charge conjugation is represented by complex conjugation of components. Notice that, as shown in \cite{bar:2007}, the KO-dimension 6 is completely natural in the Lorentzian framework. Introducing such new particles does not spoil the anomalies cancellation of the standard model. However it is not clear that this is consistent with the usual formulation of the $U(1)\times SU(2)$-symmetry.\\

A second point to understand is the formulation of the first-order condition. Indeed the definition of the first-order operators between bimodules over associative algebras is clear  \cite{ac:1994} as well as the noncommutative generalization of their symbols \cite{mdv-mas:1996a}. This notion of first-order operator  is fundamental in the noncommutative geometry approach to the standard model of \cite{cha-ac:1997}, \cite{cha-ac-mar:2007}. In contrast, the definition of first-order operator between Jordan (bi)modules is not straightforward in spite of the fact that one can write explicitely what it is for modules over the exceptional algebra $J^8_3$.\\

It seems that one should make some progress on these two points before tentative formulations of the dynamics. Then there are two natural ways to take and I think that both should be pushed. The first one is to try to adapt the spectral action principle to the present frame, the second one is to directly use the differential calculus and the connections over Jordan modules. The first approach, if it is possible, would probably be the most economical one. However two difficulties are that the spectral action principle is formulated in the frame of the Euclidean signature (instead of the Lorentzian one) and that there is the occurrence of the exceptional Jordan algebra here which needs some care. The second approach which is closer to the usual formulation of gauge theory will fully use the development of Sections 4, 5 and 6.\\

This work is currently in progress. In any case it is an occasion for starting to develop the differential calculus over Jordan algebras.

 \section*{Acknowedgements}
 
 From the beginning of my interest in the subject, Raymond Stora gave me useful advices and related references.\\
 I thank Jean Iliopoulos for his kind interest and the discussions we have had for several years.\\
 I also thank Shane Farnsworth, Gianni Landi, Todor Popov and Ivan Todorov for discussions.

%

\end{document}